\documentclass[12pt]{article}
\usepackage{amsmath}

\usepackage{amssymb}
\usepackage{mathrsfs}
\makeatletter
 
 \@addtoreset{equation}{section}
\makeatother

\newcommand{\prf}{\noindent{\bf Proof. }}

\newcommand{\rem}{\noindent{\bf Remark. }}

\newcommand{\ff}[1]{%
    {\bf F}_{#1}}

\newcommand{\qed}{\hbox{\rule[-2pt]{5pt}{11pt}}}

\newtheorem{dfn}{Definition}[section]
\newtheorem{thm}[dfn]{Theorem}

\newtheorem{lem}[dfn]{Lemma}

\newtheorem{conj}[dfn]{Conjecture}

\begin{document}
\title{On the Riemann hypothesis for self-dual weight enumerators of genus three}
\author{Koji Chinen\footnotemark[1] and Yuki Imamura\footnotemark[7]}
\date{}
\maketitle

\begin{abstract}
In this note, we give an equivalent condition for a self-dual weight enumerator of genus three to satisfy the Riemann hypothesis. We also observe the truth and falsehood of the Riemann hypothesis for some families of invariant polynomials. 
\end{abstract}

\footnotetext[1]{Department of Mathematics, School of Science and Engineering, Kindai University. 3-4-1, Kowakae, Higashi-Osaka, 577-8502 Japan. E-mail: chinen@math.kindai.ac.jp}
\footnotetext[7]{interprism Inc. Kurihara BLD 2F Nihonbashi-kakigara-cho 2-12-8 Chuo-ku Tokyo, 103-0014 Japan. E-mail:  yuki-i-xyz@outlook.jp}


\section{Introduction}\label{section:intro}
Zeta functions for linear codes were introduced by Iwan Duursma \cite{Du1} in 1999 and they have attracted attention of many mathematicians: 
\begin{dfn}\label{dfn:zeta}
Let $C$ be an $[n,k,d]$-code over $\ff q$ ($q=p^r$, $p$ is a prime) with the Hamming weight enumerator $W_C(x,y)$. Then there exists a unique polynomial $P(T)\in{\bf R}[T]$ of degree at most $n-d$ such that
\begin{equation}\label{eq:zeta_duursma}
\frac{P(T)}{(1-T)(1-qT)}(y(1-T)+xT)^n=\cdots +\frac{W_C(x,y)-x^n}{q-1}T^{n-d}+ \cdots.
\end{equation}
We call $P(T)$ and $Z(T)=P(T)/(1-T)(1-qT)$ the zeta polynomial and the zeta function of $W(x,y)$, respectively. 
\end{dfn}
If $C$ is self-dual, then $P(T)$ satisfies the functional equation (see \cite[\S 2]{Du2}): 
\begin{thm}\label{thm:func_eq}If $C$ is self-dual, then we have
\begin{equation}\label{eq:func_eq}
P(T)=P\Bigl(\frac{1}{qT}\Bigr)q^g T^{2g},
\end{equation}
where $g=n/2+1-d$. 
\end{thm}
The number $g$ is called the {\it genus} of $C$. It is appropriate to formulate the Riemann hypothesis as follows: 
\begin{dfn}\label{dfn:RH}
The code $C$ satisfies the Riemann hypothesis if all the zeros of $P(T)$ have the same absolute value $1/\sqrt q$. 
\end{dfn}
The reader is referred to \cite{Du3} and \cite{Du4} for other results by Duursma. 

\medskip
\rem The definition of the zeta function can be extended to much wider classes of invariant polynomials: let $W(x,y)$ be a polynomial of the form
\begin{equation}\label{eq:fwe_form}
W(x,y)=x^n+\sum_{i=d}^n A_i x^{n-i} y^i\in {\bf C}[x,y]\quad (A_d\ne 0)
\end{equation}
which satisfy $W^{\sigma_q}(x,y)=\pm W(x,y)$ for some $q\in {\bf R}$, $q> 0$, $q\ne1$, where
\begin{equation}\label{eq:macwilliams}
\sigma_q=\frac{1}{\sqrt{q}}\left(\begin{array}{rr} 1 & q-1 \\ 1 & -1 \end{array}\right)\quad(\mbox{the MacWilliams transform})
\end{equation}
and the action of $\sigma=\left(\begin{array}{cc} a & b \\ c & d \end{array}\right)$ on a polynomial $f(x,y)\in {\bf C}[x,y]$ is defined by $f^\sigma(x,y)=f(ax+by, cx+dy)$. Then we can formulate the zeta function and the Riemann hypothesis for $W(x,y)$ in the same way as Definitions \ref{dfn:zeta} and \ref{dfn:RH}. For the results in this direction, the reader is referred to \cite{Ch1}--\cite{Ch5}, for example. We should also note that we must assume $d, d^\perp\geq 2$, where $d^\perp$ is defined by $W^{\sigma_q}(x,y)=B_0 x^n+B_{d^\perp} x^{n-d^\perp}y^{d^\perp}+\cdots$, when considering the zeta function of $W(x,y)$. 

\medskip
We do not know much about the Riemann hypothesis for self-dual weight enumerators, but one of the remarkable results is the following theorem by Nishimura \cite[Theorem 1]{Ni}, an equivalent condition for a self-dual weight enumerator of genus one to satisfy the Riemann hypothesis: 
\begin{thm}[Nishimura]\label{thm:nishimura_g=1}
A self-dual weight enumerator $W(x,y)=x^{2d}+A_d x^d y^d+\cdots$ satisfies the Riemann hypothesis if and only if
\begin{equation}\label{eq:nishimura_g=1}
\frac{\sqrt{q}-1}{\sqrt{q}+1}{{2d}\choose{d}}\leq A_d \leq \frac{\sqrt{q}+1}{\sqrt{q}-1}{{2d}\choose{d}}.
\end{equation}
\end{thm}
Nishimura also deduces the following, the case of genus two (\cite[Theorem 2]{Ni}): 
\begin{thm}[Nishimura]\label{thm:nishimura_g=2}
A self-dual weight enumerator $W(x,y)=x^{2d+2}+A_d x^{2d+2} y^d+\cdots$ satisfies the Riemann hypothesis if and only if the both roots of the quadratic polynomial
\begin{equation}\label{eq:nishimura_g=2_eq}
A_d X^2-\left( (d-q)A_d+\frac{d+1}{d+2}A_{d+1} \right)X 
- (d+1)(q+1)\left( A_d+\frac{A_{d+1}}{d+2} \right)+(q-1){{2d+2}\choose{d}}
\end{equation}
are contained in $[-2\sqrt{q}, 2\sqrt{q}]$. 
\end{thm}
The purpose of this article is to establish an analogous equivalent condition for the case of genus three. Our main result is the following: 
\begin{thm}\label{thm:main}
A self-dual weight enumerator $W(x,y)=x^{2d+4}+A_d x^{d+4} y^d+\cdots$ satisfies the Riemann hypothesis if and only if all the roots of the cubic polynomial
\begin{equation}\label{eq:main_g=3_eq}
f_3 X^3+f_2 X^2+f_1 X+f_0 
\end{equation}
are contained in $[-2\sqrt{q}, 2\sqrt{q}]$, where $f_i$ is defined as follows.
\begin{align*}
f_3&= A_d, \\
f_2&=(q-d)A_d-\frac{d+1}{d+4}A_{d+1},\\
f_1&=\frac{1}{2}(d^2 -2q d+d-6q)A_d+(d-q+1)\frac{d+1}{d+4}A_{d+1}+\frac{(d+1)(d+2)}{(d+3)(d+4)}A_{d+2},\\
f_0&=\frac{1}{2}(q+1)(d^2+3 d-4q+2)A_d+(q+1)(d+1)(d+2)\frac{A_{d+1}}{d+4}\\
  &\quad +(q+1)\frac{(d+1)(d+2)}{(d+3)(d+4)}A_{d+2}-(q-1)\binom{2d+4}{d+4}.
\end{align*}
\end{thm}
By this theorem, we can verify the truth of the Riemann hypothesis of $W(x,y)$ only by three parameters $A_d, A_{d+1}, A_{d+2}$ (the number of parameters which are needed coincides with the genus $g$, see \cite{Ni}). Moreover, in many cases, we have $A_{d+1}=0$ and the verification of the Riemann hypothesis is simplified. 

Theorem \ref{thm:main} leads us to the consideration of the truth or falsehood of the Riemann hypothesis as the numbers $q$ and $n$ vary. As was mentioned in Remark before, $q$ can take other numbers than prime powers. In this context, we can notice the tendency that the Riemann hypothesis becomes harder to hold if $n$ or $q$ are larger. Some of the results in \cite{Ch3} and \cite{Ch4} also support it. Theorem \ref{thm:main} can illustrate this tendency by considering a certain sequence of invariant polynomials, that is 
\begin{equation}\label{eq:sequence}
W_{n,q}(x,y)=(x^2+(q-1)y^2)^n.
\end{equation}

In Section 2, we give a proof of Theorem \ref{thm:main}. In Section 3, we observe the behavior of $W_{n,q}(x,y)$, give some theoretical and experimental results, and state a conjecture on their Riemann hypothesis. 
\section{Proof of Theorem \ref{thm:main}}\label{section:proof_main}

Let $W(x,y)=x^{2d+4}+\sum_{i=d}^{2d+4}A_i x^{2d+4-i}y^i$ be a self-dual weight enumerator. Using the functional equation (\ref{eq:func_eq}) (note that $g=3$ in our case), we can assume that the zeta polynomial $P(T)$ of $W(x,y)$ is of the form
\[
P(T)= a_0 +a_1 T +a_2 T^2+a_3 T^3 +a_2 q T^4+ a_1 q^2 T^5+ a_0 q^3 T^6.
\]
We obtain another expression of $P(T)$ because $1/q\alpha$ is a root of $P(T)$ if $P(\alpha)=0$:
\begin{equation} \label{eq:P(T)factor}
P(T)=a_0 q^3 \prod_{i=1}^3 (T^2+b_i T+1/q).
\end{equation}
Comparing the coefficients, we get
\begin{align*}
b_1 +b_2 + b_3 &=a_1/a_0 q,\\
b_1 b_2 +b_2 b_3 + b_3 b_1 &=(a_2 - 3 a_0 q)/a_0 q^2,\\
b_1 b_2 b_3 &= (a_3 -2 a_1 q)/a_0 q^3.
\end{align*}
Thus $b_i$ is the roots of the following cubic polynomial:
\begin{equation} \label{eq:originaleq}
a_0 q^3 X^3 -a_1 q^2 X^2 +(a_2 -3 a_0 q)q X -a_3+2 a_1 q. 
\end{equation}
Considering the distribution of the roots of each factor $T^2+b_i T+1/q$ in (\ref{eq:P(T)factor}), we can see that a self-dual weight enumerator $W(x,y)$ of genus three satisfies the Riemann hypothesis if and only if $b_1, b_2$ and $  b_3$ are contained in $[-2/\sqrt{q}, 2/\sqrt{q}]$. By change of variable in (\ref{eq:originaleq}), we get the following: 
\begin{lem}\label{lem:equiv}
$W(x,y)$ satisfies the Riemann hypothesis if and only if all the roots of the polynomial
\begin{equation} \label{eq:originaleq-replaced}
a_0 X^3 -a_1 X^2 +(a_2 -3 a_0 q) X -a_3+2 a_1 q 
\end{equation}
are contained in $[-2\sqrt{q}, 2\sqrt{q}]$. 
\end{lem}
Our next task is to express the coefficients $a_i$ in (\ref{eq:originaleq-replaced}) by way of $A_i$ in $W(x,y)$. This can be done by comparing the coefficients of the both sides in (\ref{eq:zeta_duursma}). Our method is similar to that of Nishimura \cite{Ni}. The result is the following (here, $\alpha_{d+i} = A_{d+i}/(q-1) \binom{n}{d+i}$):
\begin{align*}
&a_0 =\alpha_{d},\\
&a_1=(d-q)\alpha_{d}+\alpha_{d+1},\\
&a_2=\frac{1}{2}d(d-2 q+1)\alpha_{d}+(d-q+1)\alpha_{d+1}+\alpha_{d+2},\\
&a_3=\frac{1}{6}d(d+1)(d-3q+2)\alpha_{d}+\frac{1}{2}(d+1)(d-2q+2)\alpha_{d+1}\\
&\qquad+(d-q+2)\alpha_{d+2}+\alpha_{d+3}.
\end{align*}
The coefficient $a_3$ is expressed by four parameters $A_d, \cdots, A_{d+3}$. By invoking the binomial moment, the number of parameters is reduced to three. In fact, we have the following: 
\begin{lem}\label{lem:moment}
Let $W(x,y)$ be a self-dual weight enumerator of the form (\ref{eq:fwe_form}) and we assume $g=3$. Then we have
\begin{equation}\label{eq:moment_g3}
\sum_{i=d+1}^{d+3} A_{i}\binom{2d+4-i}{d+1}=q \sum_{i=0}^{d+1}A_{i}\binom{2d+4-i}{d+3}.
\end{equation}
\end{lem}
\prf The equalities satisfied by the binomial moment of $W(x,y)$ is given by 
\begin{equation}\label{eq:moment}
\sum_{i=0}^{n-j} {{n-i}\choose{j}} A_i = q^{\frac{n}{2}-j} \sum_{i=0}^j {{n-i}\choose{n-j}} A_i
\qquad(j=0,1,\cdots, n)
\end{equation}
(see \cite[p.131, Problem (6)]{MaSl}). We get (\ref{eq:moment_g3}) by putting $n=2 d+4$ and $j=d+1$. \qed

\medskip
\noindent Using $A_0=1$, $A_1=\cdots =A_{d-1}=0$, we can see that (\ref{eq:moment_g3}) gives a linear relation among $A_d, \cdots , A_{d+3}$, so we can express $A_{d+3}$ by $A_d$, $A_{d+1}$ and $A_{d+2}$. Thus we get
\[
a_3 = -\frac{1}{2}(d+1)(dq+d-2q+2)\alpha_{d}-(q d+d+2)\alpha_{d+1}-(q+1)\alpha_{d+2}+1.
\]
Rewriting (\ref{eq:originaleq-replaced}) using above $a_i$, we obtain Theorem \ref{thm:main}. \qed

\section{Some examples and observations}\label{section:examples}
We examine the polynomials (\ref{eq:sequence}), which has essentially only one parameter $q$ and is easy to see the phenomenon. Using Nishimura's results ($g=1,2$) and our theorem ($g=3$), we can see that the range of $q$ for which the Riemann hypothesis is true are the following: 
\begin{eqnarray*}
g=1: & & 4-2\sqrt{3}\ (\approx 0.53590) \leq q \leq 4+2\sqrt{3}\ (\approx 7.46410) \quad(q\ne1),\\
g=2: & & -4+2\sqrt{5}\ (\approx 0.47214) \leq q \leq \alpha^2\ (\approx 3.46812)  \quad(q\ne1),
\end{eqnarray*}
where
$$\alpha=\frac{1}{6}\left(1+\sqrt[3]{5(29+6\sqrt{6})}+\sqrt[3]{5(29-6\sqrt{6})}\right),$$
and
\begin{equation}\label{eq:estimate_g=3}
g=3: \qquad \beta_1\ (\approx 0.47448) \leq q \leq \beta_3^2\ (\approx 2.47607)  \quad(q\ne1),
\end{equation}
where $\beta_1$ is the unique real root of the polynomial 
$$100t^5+495t^4+2056t^3-2928t^2+1408t-256$$
and $\beta_3$ is the positive root of the polynomial
$$13t^4+4t^3-20t^2-24t-8.$$
The cases $g=1$ and $2$ are not very complicated, but the last case needs some explanation. The relevant coefficients of $W_{4,q}(x,y)$ are 
$$A_d=A_2=4(q-1),\ A_3=0,\ A_4=6(q-1)^2.$$
Using these values, we get the explicit form of the polynomial (\ref{eq:main_g=3_eq}) as follows:
\begin{equation}\label{eq:g=3explicit}
g(X):=5X^3+5(q-2)X^2-2(11q-6)X-7q^2+20q-8.
\end{equation}
Let $D_g$ be the discriminant of $g(X)$, $X_1$ and $X_2$ be the roots of $g'(X)$ (we assume $X_1,X_2$ are real and $X_1\leq X_2$). Then, by Theorem \ref{thm:main}, $W_{4,q}(x,y)$ satisfies the Riemann hypothesis if and only if
\begin{eqnarray*}
& & D_g\geq0,\\
& & -2\sqrt{q}\leq X_1,\ X_2\leq 2\sqrt{q},\\
& & g(-2\sqrt{q})\leq 0,\ g(2\sqrt{q})\geq 0.
\end{eqnarray*}
We have 
$$\frac{D_g}{35}=100q^5+495q^4+2056q^3-2928q^2+1408q-256,$$
so $D_g\geq0$ is equivalent to
\begin{equation}\label{eq:cond_q_1}
q\geq \beta_1
\end{equation}
with the above mentioned $\beta_1$. The roots $X_i$ are given by 
\begin{eqnarray*}
X_1&=& \frac{-5(q-2)-\sqrt{25q^2+230q-80}}{15},\\
X_2&=& \frac{-5(q-2)+\sqrt{25q^2+230q-80}}{15}.
\end{eqnarray*}
The range of $q$ satisfying $-2\sqrt{q}\leq X_1$ is (note that we also have $25q^2+230q-80\geq0$)
\begin{equation}\label{eq:cond_q_2}
\frac{\sqrt{609}-23}{5}\leq q \leq \beta_2,
\end{equation}
where $\beta_2$ is the square of the unique real root of the polynomial
\begin{equation}\label{eq:cubic_beta2}
10t^3-19t^2-20t-6
\end{equation}
(this polynomial comes from the equation $-2\sqrt{q}= X_1$). The explicit value is 
$$\beta_2=\frac{1}{300}\left(761 +\sqrt[3]{386669681+396000\sqrt{17318}} +\sqrt[3]{386669681-396000\sqrt{17318}}\right)$$
($\beta_2\approx 7.38366$, this expression of $\beta_2$ can be obtained by constructing the cubic polynomial having the squares of roots of (\ref{eq:cubic_beta2}) as its roots: $100t^3-761t^2+172t-36$). The inequality $X_2\leq 2\sqrt{q}$ gives $(\sqrt{609}-23)/5\leq q$. Finally, putting $\sqrt{q}=t$, we have
\begin{eqnarray*}
g(-2\sqrt{q})&=&13t^4+4t^3-20t^2-24t-8,\\
g(2\sqrt{q})&=&13t^4-4t^3-20t^2+24t-8.
\end{eqnarray*}
The inequalities $g(-2\sqrt{q})\leq 0$ and $g(2\sqrt{q})\geq 0$ give
\begin{equation}\label{eq:cond_q_3}
0\leq q \leq \beta_3^2 \quad \mbox{and}\quad q\geq \beta_4^2 \approx 0.356397,
\end{equation}
respectively. Gathering the inequalities (\ref{eq:cond_q_1}), (\ref{eq:cond_q_2}) and (\ref{eq:cond_q_3}), we obtain the estimate (\ref{eq:estimate_g=3}). 

We can see from the above estimation that the range of $q$ for which the Riemann hypothesis is true becomes smaller as $n$ becomes larger. We show some results of numerical experiment for $W_{n,q}(x,y)$. In the following table, ``RH true'' means the range of $n$ where the Riemann hypothesis for $W_{n,q}(x,y)$ seems to be true:

\begin{center}
\begin{tabular}{c|c}
$q$ & RH true \\
\hline
2 & $2\leq q \leq 6$ \\[2mm]
$\frac{3}{2}$ & $2\leq n \leq 8$ \\[2mm]
$\frac{11}{10}$ & $2\leq n \leq 36$ \\[2mm]
$\frac{21}{20}$ & $2\leq n \leq 71$ \\[2mm]
\hline\\[-3mm]
$\frac{4}{5}$ & $2\leq n \leq 29$ \\[2mm]
$\frac{1}{2}$ & $2\leq n \leq 5$ \\[2mm]
\hline
\end{tabular}
\end{center}
These numerical examples also support the above observation. We conclude the manuscript with the following conjecture: 
\begin{conj}\label{conj:RH}
For any $n\geq 2$, there exists $q$ ($q\approx 1$) and $W_{n,q}(x,y)$ satisfies the Riemann hypothesis. 
\end{conj}

\end{document}